\documentclass[11pt, epsf]{amsart}
\usepackage{epsfig}
\usepackage{eucal}
\usepackage{subfig}
\usepackage{amssymb}
\usepackage{latexsym}
\usepackage{amsfonts}
\usepackage[latin1]{inputenc}
\usepackage{slashbox}
\usepackage{colortbl}
\usepackage{mathrsfs}

\newenvironment{namelist}[1]{%
\begin{list}{}
 {
   
   \settowidth{\labelwidth}{#1}
   \setlength{\leftmargin}{1.1\labelwidth}
  }
 }{%
\end{list}}

\newcommand{\bi}{\begin{namelist}}
\newcommand{\dis}{\displaystyle}
\newcommand{\ei}{\end{namelist}}

\newcommand{\eps}{\varepsilon}

\newtheorem{Theo}{Theorem}

\title[FEM For Solving The Dirac Eigenvalue Problem With Linear Basis Functions]{Finite Element Method For Solving The Dirac Eigenvalue Problem With Linear Basis Functions}
\author{Hasan Almanasreh}
\thanks{Mathematics Department, Hebron University, P.O. Box 40, Hebron, West Bank, Palestine}

\begin{document}

\maketitle

\begin{abstract} In this work we will treat the spurious eigenvalues obstacle that appears in the computation of the radial Dirac eigenvalue problem using numerical methods. The treatment of the spurious solution is based on applying Petrov-Galerkin finite element method. The significance of this work is the employment of just continuous basis functions, thus the need of a continuous function which has a continuous first derivative as a basis, as in \cite{ALMAN,ALMANA,ALMA}, is no longer required. The Petrov-Galerkin finite element method for the Dirac eigenvalue problem strongly depends on a stability parameter, $\tau$, that controls the size of the diffusion terms added to the finite element formulation for the problem. The mesh-dependent parameter $\tau$ is derived based on the given problem with the particular basis functions.\\

\end{abstract}

\section{Introduction}

In quantum mechanics, the Dirac partial differential equation describes the relativistic behavior of the electrons around the nucleus. That is, the energies (eigenvalues) of the electron in the orbital levels can be computed by solving the Dirac eigenvalue problem. The eigenvalues of an electron in the many-electron systems (nucleus with more than one electron) can be approximated based on the single-electron systems (nucleus with just one electron around), see \cite{LIND,MOH,ROS,SHAB1}. The main obstacle in solving the Dirac eigenvalue problem is that the genuine eigenvalues are polluted by spurious eigenvalues (called spectrum pollution) \cite{ACK,ALMAN,ALMA,PES,TUPS}. The spurious solutions appeared in many numerical computations of eigenvalue problems \cite{DAW,MUR,SCH,ZHA}. On the other hand, applying the numerical methods, with their general forms, to the Dirac eigenvalue problem encountered the presence of spurious eigenvalues; such these numerical methods are B-spline method \cite{FISCHP,FISCHZ,JOH,SHAB}, FEM \cite{ALMANA,ALMA,MUL}, meshfree method \cite{ALMAN}, and  FDM \cite{SAL}.
In this work, we present a stable computation, using the FEM, of the eigenvalues of the Dirac operator by means of a complete remedy of the spectrum pollution. To present the work, consider first the free Dirac operator $\mathbf{H}_0$ with the Coulomb potential $V$
\begin{equation}\label{1}
\mathbf{H}=\mathbf{H}_0+V\, .
\end{equation}
The free operator $\mathbf{H}_0:H^1({\mathbb R}^3;{\mathbb C}^4)\to L^2({\mathbb R}^3;{\mathbb C}^4)$ is given by
\begin{equation}\label{2}
\mathbf{H}_0 = -i{\hslash}c\boldsymbol{\alpha}\cdot \boldsymbol{\nabla} + mc^2 \beta\, ,
\end{equation}
where $\hslash$ is the Planck constant divided by $2\pi$, the operator $\boldsymbol{\nabla}\!=\!(\frac{\partial}{\partial x_1},\frac{\partial}{\partial x_2},\frac{\partial}{\partial x_3})$, $c$ is the speed of light, and $m$ is the electron mass at rest, the symbols $\boldsymbol{\alpha}\!=\!(\alpha_1,\alpha_2,\alpha_3)$ and $\beta$ are the $4\times 4$ Dirac matrices given by
$$
\alpha_j = \left(
\begin{array}{cc}
0 & \sigma_j \\
\sigma_j & 0
\end{array}
\right)\;\;\text{and}\;\;
\beta = \left(
\begin{array}{cc}
I & 0 \\
0 & -I
\end{array}
\right)\, .
$$
Here $I$ and $0$ are the $2\times 2$ identity and zeros matrices respectively, and $\sigma_j$'s are the $2\times 2$ Pauli matrices
\vspace{-0.2cm}
$$
\sigma_1 = \left(
\begin{array}{cc}
0 & 1 \\
1 & 0
\end{array}
\right),\;\;
\sigma_2 = \left(
\begin{array}{cc}
0 & -i \\
i & 0
\end{array}
\right)\, ,
\;\;\text{and}\;\;
\sigma_3 = \left(
\begin{array}{cc}
1 & 0 \\
0 & -1
\end{array}
\right)\, .
$$
The Coulomb potential $V$ is a multiplicative operator given by
\begin{equation}\label{3}
 V(x)\!=\!\frac{-z}{|x|}I,
\end{equation}
here $I$ is the $4\times 4$ identity matrix, where $I$ will be dropped for simplicity. The independent variable $x=(x_1,x_2,x_3)$ and $z\in\{1,2,\ldots 137\}$ is the electric charge number. The operator $\mathbf{H}_0$ is essentially self-adjoint on $C_0^\infty(\mathbb{R}^3;\mathbb{C}^4)$ and self-adjoint on $H^1({\mathbb R}^3;{\mathbb C}^4)$. Thus, the whole operator $\mathbf{H}$ is self-adjoint on $H^1({\mathbb R}^3;{\mathbb C}^4)$. Moreover, the spectrum of $\mathbf{H}$ is $(-\infty,-mc^2]\cup\{\lambda^k\}_{k\in\mathbb N}\cup[mc^2,+\infty)$, where $\{\lambda^k\}_{k\in\mathbb N}$ is a discrete sequence of eigenvalues (relativistic energies).

The Dirac eigenvalue problem is given by
\begin{equation}\label{4}
\mathbf{H}u(x)=\lambda u(x)\, ,
\end{equation}
where $u\in H^1({\mathbb R}^3;{\mathbb C}^4)$. Usually, the radial Dirac operator is considered when the computation of the eigenvalues $\lambda$ is concerned. The radial operator can be obtained by separation of variables of the radial and angular parts. That is, by assuming
$u(x)=\displaystyle\frac{1}{r}\left(
\begin{array}{c}
f(r)\mathscr{Z}_{\kappa,m}(\varpi,\theta) \\
i\,g(r)\mathscr{Z}_{-\kappa,m}(\varpi,\theta)
\end{array}
\right)$, where $r=|x|$ is the radial variable, $f$ and $g$ are the Dirac large and small radial functions respectively, $\mathscr{Z}_{\cdot,m}$ is the angular part of the wave function $u$, and $\kappa$ is the spin-orbit coupling parameter defined as $\kappa\!=\!(-1)^{\jmath+\ell+\frac{1}{2}}(\jmath+\frac{1}{2})$, where $\jmath$ and $\ell$ are the total and orbital angular momentum numbers respectively. By this separation, the radial Dirac eigenvalue problem is then given by, see, e.g., \cite{THA},
\begin{equation}\label{5}
H_\kappa\varphi(r)=\lambda\varphi(r)\, ,\quad \text{where}
\end{equation}
\begin{equation}\label{6}
H_\kappa = \dis\left(\dis
\begin{array}{cc}
\dis mc^2+V(r) &\dis c\big(\!-\!\frac{d}{dr}+\frac{\kappa}{r}\big) \\
\dis c\big(\frac{d}{dr}+\frac{\kappa}{r}\big) &\dis -mc^2+V(r)
\end{array}
\right)\;\; \text{and}\;\; \varphi(r) = \left(
\begin{array}{c}
f(r) \\
g(r)
\end{array}
\right)\, .
\end{equation}
As defined before, $\lambda$ is the relativistic energy, and $V(r)\!=\!-z/r$ is the radial Coulomb potential.\\

The radial Dirac operator is a convection-dominated operator, see, e.g., \cite{ALM,DES,IDE,LIN}, which causes instability in the numerical approximation of the eigenvalues. That is, the presence of the gradient in the off diagonal of the operator $H_\kappa$ and the absence of the Laplace operator is the core of the spuriosity problem in the numerical computation \cite{ALMAN,ALMA}.\\

In this work we will provide a stable finite element computation of the eigenvalues $\lambda$ of the operator $H_\kappa$. The finite element scheme we provide here based on applying the stream line upwind Petrov-Galerkin (SUPG) instead of the usual Galerkin FEM to produce diffusivity, controlled by a stability parameter ($\tau$) derived for the specific problem. The parameter $\tau$ controls the size of the added diffusion terms to the usual Galerkin formulation of the problem and its derivation is particular for finite element formulation of the Dirac eigenvalue problem with linear basis functions. The derivation of $\tau$ here is simpler than in \cite{ALMAN,ALMA}, moreover, the need of the accumulation of the eigenvalues \cite{GRI} is not required in the derivation.\\

The paper is arranged as follows; in Section 2 we provide some required preliminaries. In Section 3, we talk about the Galerkin and the SUPG finite element formulation to the problem and discuss the scheme of stability. The derivation of the stability parameter $\tau$ is treated in Section 4. Finally, we support our work by computational results in Section 5 and provide a discussion.

\section{Preliminaries}

\subsection{The function space}
In \cite{ALMAN,ALMA} additional requirement of the functional space is considered, that is, the radial Dirac functions $f$ and $g$ are assumed to be $C^1$ (the space of continuous functions that have continuous first derivatives). This requirement is time consuming in the computation, so, in this work we show that this requirement is no longer needed. To determine the specific function space, firstly, it is clear that the radial functions $f$ and $g$ belong to the space $H^1(\Omega)$, where $\Omega=[0,\infty)$ (The radial domain). Also, the functions $f$ and $g$ should vanish near the boundaries (close to and far away from the nucleus), so homogeneous Dirichlet boundary condition is considered. Therefore, $f,g\in\displaystyle H_0^1(\Omega)$. It should also be notified that $f$ and $g$ must smoothly vanish at the boundaries (in a damping way) for all states except $1s_{1/2}$ and $2p_{1/2}$ ($\kappa=-1$ and $\kappa=1$ respectively). That is, for better approximation, homogeneous Neumann boundary condition should also be applied in any computation of the eigenvalues for these states. In the presented work, general and unified treatment is considered for the boundary conditions, that is, homogeneous Dirichlet boundary condition is only assumed throughout all computations, for more readings see \cite{ALMAN,ALMA}.

\subsection{Extended nucleus}
It is notable that the Coulomb potential is singular near $r=0$, so careful treatment should be taken into account to avoid this singularity. That is, extended nucleus is considered in this case. Extended nucleus means that to assume another distribution (a function that has to be at least $C^1-$function) of the electric charge on the domain $[0,R]$ ($R$ is the radius of the nucleus) while keeping the Coulomb potential on the rest of the domain $[R,\infty)$. The distribution of the electric charge on $[0,R]$ can be, e.g., Fermi or uniform distribution see \cite{ALMAN,ALMA}. In this work, we will assume a uniformly distributed charge along the interval $[0,R]$. Computationally, we first treat point nucleus where a cut-off domain, $[R,\infty)$, is considered to avoid the singularity, thenafter we extend the computation on the whole domain $[0,\infty)$. For the point nucleus case, we can test both the convergence of the genuine eigenvalues and the remedy of the spectrum pollution. This is because for point nucleus we can compare our results with the exact values of the eigenvalues that are given by the relativistic formula
\begin{equation}\label{7}
\lambda_{n_r,\kappa}=\frac{mc^2}{\sqrt{1+\frac{z^2\gamma^2} {(n_r-1+\sqrt{\kappa^2-z^2\gamma^2})^2}}}\,,
\end{equation}
where $\gamma$ is the fine structure constant which has the value $1/c$ in atomic unit, and $n_r=1,2,\ldots$ is the orbital level number. To make the comparison simpler, the exact eigenvalues $\lambda_{n_r,\kappa}$ and the computed ones are shifted by $-mc^2$.

\subsection{Exponentially distributed nodes}
Since the wave functions oscillates heavily close to the nucleus compared to the regions away from it, more data is required at this region to get more accurate approximation. For this reason, exponential distribution of the nodes is considered. Here, the nodes are distributed along the interval of computation $[a, b]$ by the following formula
\begin{equation}\label{8}
r_i=\displaystyle\exp\Big({\ln(a+\eps)+\big(\frac{\ln(b+\eps)- \ln(a+\eps)}{n+1}\big)i}\Big)-\eps\, , \;\;\; i=0,1,2,\ldots,n+1,
\end{equation}
where $n+1$ is the number of subintervals and $\eps\in[0\,,\,1]$ is the nodes intensity parameter \cite{ALMAN}. The goal of introducing the parameter $\eps$ is to control the intensity of the nodal points near the nucleus. As smaller $\eps$ as more nodes dragged closed to the nucleus and vice versa. In \cite{ALMAN} a study is performed about the best choices of $\eps$, where it is shown that the most appropriate values of $\eps$ are those that are living in the interval $[10^{-6}\,,\,10^{-4}]$.

\section{The Petrov-Galekin formulation}
Recall the radial Dirac eigenvalue problem; find $(\lambda,\,(f,g))\in \mathbb{R}\times(H_0^1([0,\infty]))^2$ such that
\begin{equation}\label{9}
H_\kappa\varphi(r)=\lambda\varphi(r)\, ,\quad \text{where}
\end{equation}
\begin{equation}\label{10}
H_\kappa = \dis\left(\dis
\begin{array}{cc}
\dis mc^2+V(r) &\dis c\big(\!-\!\frac{d}{dr}+\frac{\kappa}{r}\big) \\
\dis c\big(\frac{d}{dr}+\frac{\kappa}{r}\big) &\dis -mc^2+V(r)
\end{array}
\right)\;\; \text{and}\;\; \varphi(r) = \left(
\begin{array}{c}
f(r) \\
g(r)
\end{array}
\right)\, .
\end{equation}
To discretise the problem, let $V$ be the space of continuous functions and $V^L$ be the subspace of $V$ that consists of continuous linear polynomials. Assume a partition $K_h$ consisting of exponentially distributed points in $[a, b]$. Now, let $V_h^L$ be the finite subspace of $V^L$ consisting of piecewise continuous linear polynomials spanned by the below linear functions $\phi_j$ on the partition $K_h$
\begin{displaymath}
\phi_j(r) = \left\{ \begin{array}{ll}
\displaystyle\frac{r\,-\,r_{j-1}}{h_j}\, ,& r\in [r_{j-1}, r_j]\, , \\
\displaystyle\frac{r_{j+1}\,-\,r}{h_{j+1}}\, ,& r\in [r_j, r_{j+1}]\, ,\\
0\, ,&\text{elsewhere},
\end{array} \right.
\end{displaymath}
where $j=0,1,2,\ldots,n+1$. Now, if $f, g\in V_h^L$, then
\begin{equation}\label{11}
f(r)=\sum_{j=0}^{n+1}\zeta_j\phi_j(r)\, ,
\end{equation}
\vspace{-0.3cm}
\begin{equation}\label{11half}
g(r)=\sum_{j=0}^{n+1}\xi_j\phi_j(r)\, ,
\end{equation}
where $\zeta_j$ and $\xi_j$ are respectively the values of the functions $f$ and $g$ at the node $r_j$. To construct the Galerkin FEM for the problem, we assume that $f,g\in V_h^L$. Since homogeneous boundary condition is assumed, then $f$ and $g$ should vanish at the boundaries, that is $\zeta_0=\zeta_{n+1}=\xi_0=\xi_{n+1}=0$. Now the Galerkin FEM is read as to multiply (\ref{9}) by test functions $(v, 0)^t$ and $(0, v)^t$ and integrate over the whole domain $\Omega$ this gives the weak form of the problem
\begin{align}\label{12}
\sum_{j=1}^n\Big(\dis\int_\Omega(mc^2+V(r))\phi_j(r) v(r) dr\Big)\zeta_j+\sum_{j=1}^n\Big(\dis\int_\Omega(-c\phi_j'(r)+\frac{c\kappa}{r}\phi_j(r)) v(r)dr\Big)\xi_j\\\nonumber
=\lambda\sum_{j=1}^n\dis\Big(\int_\Omega\phi_j(r)v(r) dr\Big)\zeta_j,
\end{align}
and
\begin{align}\label{13}
\sum_{j=1}^n\Big(\dis \int_\Omega (c\phi_j'(r)+\frac{c\kappa}{r}\phi_j(r)) v(r) dr\Big)\zeta_j +\sum_{j=1}^n\Big(\dis \int_\Omega (-mc^2+V(r))\phi_j(r) v(r)dr\Big)\xi_j\\\nonumber
=\lambda\sum_{j=1}^n\Big(\dis\int_\Omega\phi_j(r) v(r)dr\Big)\xi_j.
\end{align}
To complete the numerical formulation, let the test function be an element of the same space $\mathcal{V}^l_h$ such that $v=\phi_i$, $i=1,2,\ldots,n$, in (\ref{12}) and (\ref{13}), this provides
\begin{equation}\label{14}
A\,Y=\lambda B\,Y\, .
\end{equation}
Equation (\ref{14})is a symmetric generalized eigenvalue problem, where $A$ and $B$ are both $2n\times 2n$ symmetric block matrices defined as
\begin{equation}\label{15}
A=\left[\begin{array}{c|c}
mc^2\dis[a_{ij}]_{000}+\dis[a_{ij}]_{000}^V & -c\dis[a_{ij}]_{010}+c\kappa \dis[a_{ij}]_{001} \\
\hline
c\dis[a_{ij}]_{010}+c\kappa \dis[a_{ij}]_{001} & -mc^2\dis[a_{ij}]_{000}+\dis[a_{ij}]_{000}^V
\end{array}\right],
\end{equation}
and
\begin{equation}\label{16}
B=\left[\begin{array}{c|c}
\dis[a_{ij}]_{000}& 0 \\
\hline
0 & \dis[a_{ij}]_{000}
\end{array}\right],
\end{equation}
where $\dis[a_{ij}]_{\rho\sigma\nu}^q$ is an $n\times n$ matrix defined as
\begin{equation}\label{17}
(\dis[a_{ij}]_{\rho \sigma \nu}^q)_{ij}=\int_{\Omega}\phi_j^{(\sigma)}\,\phi_i^{(\rho)}\,r^{-\nu}\,q(r)\,dr\; ,\;\;\text{where}\;\;\phi^{(\rho)}(r)=\frac{d^\rho}{dr^\rho}\phi (r)\, .
\end{equation}

To stabilize the computation, that is to get red of the spectrum pollution, the streamline upwind Petrov-Galerkin (SUPG) FEM is considered instead of the Galerkin FEM for the problem \cite{ALMAN,ALMA,ALM,DES,IDE}. The idea of using SUPG FEM is to introduce diffusion terms in the weak formulation of the problem to stabilize the numerical computation.\\

The construction of the SUPG FEM for the radial eigenvalue problem is to promote the test function to also include its first derivative, that is to multiply (\ref{9}) by $(v, \tau v')^t$ and $(\tau v', v)^t$ instead of just $(v, 0)^t$ and $(0, v)^t$. This will introduce diffusion terms of the form $\tau \dis[a_{ij}]_{110}$, where $\dis[a_{ij}]_{110}=\dis\int_{\Omega}\phi_j'(r)\,\phi_i'(r)\,dr$, on the main diagonal of the generalized matrix $\mathcal{A}$. The parameter $\tau$ is the so-called stability parameter that always depends on the generated mesh. The role of $\tau$ is to control the size of the added diffusion terms. To formulate the SUPG FEM, multiply (\ref{9}) by $(v, \tau v')^t$ and $(\tau v', v)^t$ and integrate over the domain $\Omega$

\begin{align}\label{18}
\dis\int_\Omega (mc^2+V(r))f(r)\, v(r) dr+\dis\int_\Omega(-cg'(r)+\frac{c\kappa}{r}g(r))\, v(r) dr+\,\tau\!\!\dis\int_\Omega R_2(f(r),g(r))\,  v'(r) dr \\\nonumber
=\lambda\dis\int_\Omega f(r)\, v(r) dr\, ,
\end{align}
and
\begin{align}\label{19}
\dis\int_\Omega (cf'(r)+\frac{c\kappa}{r}f(r))\, v(r) dr+\dis\int_\Omega (-mc^2+V(r))g(r)\, v(r) dr+\,\tau\!\!\dis\int_\Omega R_1(f(r),g(r))\,  v'(r) dr \\\nonumber
=\lambda\dis\int_\Omega g(r)\, v(r) dr\, ,
\end{align}
where
\begin{equation}\label{20}
R_{1}(f(r),g(r))=(mc^2+V(r)-\lambda)f(r)-cg'(r)+\frac{c\kappa}{r}g(r)\, ,
\end{equation}
and
\begin{equation}\label{21}
R_{2} (f(r),g(r))=(-mc^2+V(r)-\lambda)g(r)+cf'(r)+\frac{c\kappa}{r}f(r)\, ,
\end{equation}

To discretize the weak formulation let $v,f,g\in V_h^L$ such that $f$ and $g$ as given by (\ref{11}) and (\ref{11half}) respectively, and $v=\phi_i$, $i=1,2,\ldots,n$, this leads to the generalized eigenvalue problem
\begin{equation}
\mathcal{A}Y=\lambda\mathcal{B}Y\, ,
\end{equation}
where $\mathcal{A}$ and $\mathcal{B}$ are given by
\begin{equation}\label{22}
\mathcal{A}=\left[\begin{array}{c|c}
mc^2 \dis[a_{ij}]_{000}+ \dis[a_{ij}]_{000}^V+ & -c \dis[a_{ij}]_{010}+c\kappa  \dis[a_{ij}]_{001}+\\
+c\tau  \dis[a_{ij}]_{110}+c\tau\kappa  \dis[a_{ij}]_{101} & -mc^2\tau  \dis[a_{ij}]_{100}+\tau  \dis[a_{ij}]_{100}^V\\
\hline
c \dis[a_{ij}]_{010}+c\kappa  \dis[a_{ij}]_{001}+  & -mc^2 \dis[a_{ij}]_{000}+ \dis[a_{ij}]_{000}^V+\\
mc^2\tau \dis[a_{ij}]_{100}+\tau \dis[a_{ij}]_{100}^V & -c\tau \dis[a_{ij}]_{110}+c\tau\kappa \dis[a_{ij}]_{101}
\end{array}\right]
\end{equation}
and
\begin{equation}\label{23}
\mathcal{B}=\left[\begin{array}{c|c}
 \dis[a_{ij}]_{000} & \tau  \dis[a_{ij}]_{100} \\
\hline
\tau \dis[a_{ij}]_{100} &  \dis[a_{ij}]_{000}
\end{array}\right)\, .
\end{equation}

\section{The stability parameter $\tau$}

\begin{Theo}
Let $\zeta_{j-1}$ and $\zeta_{j+1}$ (resp. $\xi_{j-1}$ and $\xi_{j+1}$) be the values of the radial function $f$ (resp. $g$) at the nodes $r_{j-1}$ and $r_{j+1}$ respectively. Then $\zeta_{j-1}$, $\zeta_{j+1}$, $\xi_{j-1}$, and $\xi_{j+1}$ can be approximated in the vicinity of $r$ at infinity by
\begin{eqnarray*}
 && \zeta_{j-1}\approx \zeta_j+\big(-mc h_j-\frac{h_j}{c}\lambda\big)\xi_j\, . \\
 && \zeta_{j+1}\approx \zeta_j+\big(mc h_{j+1}+\frac{h_{j+1}}{c}\lambda\big)\xi_j\, . \\
  && \xi_{j-1}\approx \xi_j+\big(-mc h_j+\frac{h_j}{c}\lambda\big)\zeta_j\, . \\
   &&  \xi_{j+1}\approx \xi_j+\big(mc h_{j+1}-\frac{h_{j+1}}{c}\lambda\big)\zeta_j\, .
\end{eqnarray*}
\end{Theo}
\hspace{-4mm}\underline{\emph{Proof}}. Given the two-equation system of (\ref{9})
\begin{eqnarray}\label{24}
\big(mc^2+V(r)\big)f(r)+c\big(-g'(r)+\frac{\kappa}{r}g(r)\big)=\lambda f(r),
\end{eqnarray}
and
\begin{eqnarray}\label{25}
c\big(f'(r)+\frac{\kappa}{r}f(r)\big)+\big(-mc^2+V(r)\big)g(r)=\lambda g(r)\,.
\end{eqnarray}
As $r$ approaches infinity, the above two equations are reduced to
\begin{eqnarray}\label{26}
mc^2f(r)-c g'(r)=\lambda f(r),
\end{eqnarray}
and
\begin{eqnarray}\label{27}
c f'(r)-mc^2 g(r)=\lambda g(r)\,.
\end{eqnarray}
To obtain the desired formula for $\zeta_{j-1}$ and $\xi_{j-1}$, we use the backward difference approximation for $f'$ and $g'$ as follows

$$\Rightarrow f'|_{r_j}\approx \frac{f(r_j)-f(r_{j-1})}{r_j-r_{j-1}}= \frac{\zeta_j-\zeta_{j-1}}{h_j},$$
and
$$\Rightarrow g'|_{r_j}\approx \frac{g(r_j)-g(r_{j-1})}{r_j-r_{j-1}}= \frac{\xi_j-\xi_{j-1}}{h_j} .$$
By these approximations of the derivatives, (\ref{26}) and (\ref{27}) at the node $r_j$ can be written as
\begin{eqnarray}\label{28}
mc^2\zeta_j-c \frac{\xi_j-\xi_{j-1}}{h_j}=\lambda \zeta_j,
\end{eqnarray}
and
\begin{eqnarray}\label{29}
c \frac{\zeta_j-\zeta_{j-1}}{h_j}-mc^2 \xi_j=\lambda \xi_j\,.
\end{eqnarray}
Simplifying (\ref{28}) and (\ref{29}) gives the desired result for $\xi_{j-1}$ and $\zeta_{j-1}$.

To obtain the corresponding formulas for $\xi_{j+1}$ and $\zeta_{j+1}$, we assume (\ref{26}) and (\ref{27}) is true for $r_j$ and then use the forward difference approximation for $f'$ and $g'$ as follows

$$\Rightarrow f'|_{r_j}\approx \frac{f(r_{j+1})-f(r_{j})}{r_{j+1}-r_{j}}= \frac{\zeta_{j+1}-\zeta_{j}}{h_{j+1}},$$
and
$$\Rightarrow g'|_{r_j}\approx \frac{g(r_{j+1})-g(r_{j})}{r_{j+1}-r_{j}}= \frac{\xi_{j+1}-\xi_{j}}{h_{j+1}} .$$\hfill{$\blacksquare$}\\

\begin{Theo}
The stability parameter $\tau$ that appears in the weak formulation (\ref{18}) and (\ref{19}) has the form
\begin{equation}\label{30}
\tau:=\tau_j=\frac{1}{3}(h_{j+1}-h_j) .
\end{equation}
\end{Theo}
\hspace{-4mm}\underline{\emph{Proof}}. Consider the weak formulation of the radial Dirac equation in the vicinity of $r$ at infinity.

\begin{eqnarray}\label{31}
(mc^2-\lambda) \dis [a_{ij}]_{000}f+\tau c [a_{ij}]_{110}f-(\tau m c^2-c+\tau\lambda) [a_{ij}]_{100}\;g=0,
\end{eqnarray}
and
\begin{eqnarray}\label{32}
(\tau mc^2-c-\tau\lambda) [a_{ij}]_{100}f-\tau c [a_{ij}]_{110}\;g-(mc^2+\lambda)[a_{ij}]_{000}\;g=0,
\end{eqnarray}
where $f=(\zeta_1,\zeta_2,\ldots,\zeta_n)$ and $g=(\xi_1,\xi_2,\ldots,\xi_n)$ are respectively the nodal values of the functions $f$ and $g$. Note that we have used that $[a_{ij}]_{010}=-[a_{ij}]_{100}$ in the above formulation. Now, using the following values of the integrals,

\begin{table}[ht]
\caption{The element integrals of the matrices $[a_{ij}]_{000}$, $[a_{ij}]_{100}$, and $[a_{ij}]_{110}$.}
\centering
\begin{tabular}{@{} ||c|c|c|c|| @{} }
\hline\hline
\backslashbox{$j^{th}$ row of}{Column} & $j-1$ & $j$ & $j+1$ \\ [0.5ex]
\hline
$[a_{ij}]_{000}$ & $\dis\frac{1}{6}h_{j}$ & $\dis\frac{1}{3}(h_j+h_{j+1})$ & $\dis\frac{1}{6}h_{j+1}$ \\
\hline
$[a_{ij}]_{100}$ & $\dis\frac{1}{2}$ & $0$ &$\dis-\frac{1}{2}$\\
\hline
$[a_{ij}]_{110}$ & $\dis-\frac{1}{h_j}$ & $\dis\frac{1}{h_j}+\frac{1}{h_{j+1}}$ & $\dis-\frac{1}{h_{j+1}}$ \\
 [1ex]
\hline
\hline
\end{tabular}
\end{table}
equations (\ref{31}) and (\ref{32}) becomes
\begin{eqnarray}\label{33}
(mc^2-\lambda) \big[\dis\frac{h_j}{6}\zeta_{j-1}+\dis\frac{1}{3}(h_j+h_{j+1})\zeta_j+ \dis\frac{h_{j+1}}{6}\zeta_{j+1}\big] +\tau c \big[\dis-\frac{1}{h_j}\zeta_{j-1}+\\\nonumber
 +\dis\frac{h_j+h_{j+1}}{h_j h_{j+1}}\zeta_j\dis-\frac{1}{h_{j+1}}\zeta_{j+1}\big]-(\tau m c^2-c+\tau\lambda) \big[\frac{1}{2}\xi_{j-1}-\frac{1}{2}\xi_{j+1}\big]=0,
\end{eqnarray}
and
\begin{eqnarray}\label{34}
(\tau mc^2-c-\tau\lambda) \big[\frac{1}{2}\zeta_{j-1}-\frac{1}{2}\zeta_{j+1}\big]-\tau c \big[\dis-\frac{1}{h_j}\xi_{j-1} +\dis\frac{h_j+h_{j+1}}{h_j h_{j+1}}\xi_j+\\\nonumber
-\dis\frac{1}{h_{j+1}}\xi_{j+1}\big]-(mc^2+\lambda)\big[\dis\frac{h_j}{6}\xi_{j-1}+\dis\frac{1}{3}(h_j+h_{j+1})\xi_j+ \dis\frac{h_{j+1}}{6}\xi_{j+1}\big]=0.
\end{eqnarray}
Using Theorem 1, the above two equations can be written as
\begin{eqnarray}\label{35}
\Big[\frac{1}{6}h_jmc^2-\frac{1}{6}h_j\lambda-\frac{\tau c}{h_j}\frac{1}{3}(h_{j+1}+h_j)mc^2-\frac{1}{3}(h_{j+1}+h_j)\lambda + \tau c\frac{h_{j+1}+h_j}{h_jh_{j+1}}+\\\nonumber
+\frac{1}{6}\,h_{j+1}\,m\,c^2-\frac{1}{6}h_{j+1}\,\lambda-\frac{\tau\, c}{h_{j+1}}+(\frac{c}{2}-\frac{\tau\, m\, c^2}{2}-\frac{\tau\, \lambda}{2})(-m\,c\,h_j+\frac{\lambda\, h_j}{c})+\\\nonumber
+(-\frac{c}{2}+\frac{\tau m c^2}{2}+\frac{\tau\lambda}{2})(mch_{j+1}-\frac{\lambda h_{j+1}}{c})\Big]\zeta_j+ \Big[(\frac{1}{6}h_jmc^2-\frac{1}{6}h_j\lambda-\frac{\tau c}{h_j})\times\\\nonumber
\times(-mch_j-\frac{\lambda h_j}{c})+(\frac{1}{6}h_{j+1}m\,c^2-\frac{1}{6}h_{j+1} \lambda-\frac{\tau c}{h_{j+1}})(m\,c\,h_{j+1}+\frac{\lambda h_{j+1}}{c})+\\\nonumber
\frac{c}{2}-\frac{\tau m c^2}{2}-\frac{\tau \lambda}{2}-\frac{c}{2}+\frac{\tau m c^2}{2}+\frac{\tau \lambda}{2}\Big]\xi_j=0,
\end{eqnarray}
and
\begin{eqnarray}\label{36}
\Big[\frac{\tau m c^2}{2}-\frac{\tau \lambda}{2}-\frac{c}{2}-\frac{\tau m c^2}{2}+\frac{\tau \lambda}{2}+\frac{c}{2}+ (\frac{\tau c}{h_{j}}
-\frac{1}{6}h_{j}mc^2-\frac{1}{6}h_{j}\lambda)(-mch_{j}+\frac{\lambda h_{j}}{c})+\\\nonumber+(\frac{\tau c}{h_{j+1}}-\frac{1}{6}h_{j+1}mc^2 -\frac{1}{6}h_{j+1}\lambda)(mch_{j+1}
-\frac{\lambda h_{j+1}}{c})\Big]\zeta_j+\Big[(\frac{\tau m c^2}{2}-\frac{\tau \lambda}{2}-\frac{c}{2})\times\\\nonumber
\times(-mch_j-\frac{\lambda h_j}{c})+(-\frac{\tau m c^2}{2}+\frac{\tau\lambda}{2}
+\frac{c}{2})(mch_{j+1}+\frac{\lambda h_{j+1}}{c})+\frac{\tau c}{h_j}-\frac{1}{6}h_jmc^2+\\\nonumber
-\frac{1}{6}h_j\lambda-\tau c\frac{h_{j+1}+h_j}{h_jh_{j+1}}
-\frac{1}{3}(h_{j+1}+h_j)mc^2 -\frac{1}{3}(h_{j+1}+h_j)\lambda+\frac{\tau c}{h_{j+1}}+\\\nonumber -\frac{1}{6}h_{j+1}mc^2-\frac{1}{6}h_{j+1}\lambda\Big]\xi_j=0.
\end{eqnarray}
Assuming $m=1$, as $c\to\infty$, and after some algebraic simplifications, equations (\ref{35}) and (\ref{36}) becomes
\begin{eqnarray}\label{37}
\Big[\frac{c^2}{6}-\frac{\lambda}{6}-\frac{c^2}{3}-\frac{\lambda}{3}-\frac{\tau c}{h_jh_{j+1}}\!+\!\frac{\tau c}{h_jh_{j+1}}-\frac{c^2}{2} \!+\!\frac{\tau c^3}{2}\!+\!\frac{\tau\lambda c}{2}
+\frac{\lambda}{2}-\frac{\tau c\lambda}{2}\Big](h_j\!+\!h_{j+1})\zeta_j+\\\nonumber+ \Big[-\frac{h_j^2c^3}{6}+\frac{h_j^2c\lambda}{6}+\tau c^2-\frac{ch_j^2\lambda}{6}
+\tau\lambda+\frac{h_{j+1}^2c^3}{6}-\frac{h_{j+1}^2 c\lambda}{6}-\tau c^2+\frac{h_{j+1}^2c\lambda}{6}-\tau\lambda\Big]\xi_j=0,
\end{eqnarray}
and
\begin{eqnarray}\label{38}
\Big[-\tau c^2+\frac{h_{j}^2c^3}{6}+\frac{h_j^2c\lambda}{6}+\tau \lambda-\frac{h_j^2c\lambda}{6}+\tau c^2-\frac{h_{j+1}^2c^3}{6}
-\frac{h_{j+1}^2\lambda c}{6}
-\tau \lambda+\frac{h_{j+1}^2c\lambda}{6}\Big]\zeta_j+\\\nonumber+\Big[-\frac{\tau c^3}{2}+\frac{\tau c\lambda}{2}+\frac{c^2}{2}- \frac{\tau c\lambda}{2} +\frac{\lambda}{2}-\frac{c^2}{6}
-\frac{\lambda}{6}+\frac{\tau c}{h_jh_{j+1}}- \frac{\tau c}{h_jh_{j+1}}-\frac{c^2}{3}-\frac{\lambda}{3}\Big](h_j+h_{j+1})\xi_j=0.
\end{eqnarray}
Dividing (\ref{37}) and (\ref{38}) by $(h_j+h_{j+1})$ yields

\begin{eqnarray}\label{39}
\Big[\frac{c^2}{6}-\frac{\lambda}{6}-\frac{c^2}{3}-\frac{\lambda}{3}-\frac{c^2}{2}+\frac{\tau c^3}{2}+\frac{\lambda}{2}\Big]\zeta_j+ \Big[\frac{c^3}{6}(h_{j+1}-h_j)]\xi_j=0,
\end{eqnarray}
and
\begin{eqnarray}\label{40}
\Big[-\frac{c^3}{6}(h_{j+1}-h_j)\Big]\zeta_j+\Big[-\frac{\tau c^3}{2}+\frac{c^2}{2}+\frac{\lambda}{2}-\frac{c^2}{6}-\frac{\lambda}{6}-\frac{c^2}{3}-\frac{\lambda}{3}\Big]\xi_j=0.
\end{eqnarray}
Simplifying the above two equations provides
\begin{eqnarray}\label{41}
\Big[\frac{\tau c^3}{2}\Big]\zeta_j+ \Big[\frac{c^3}{6}(h_{j+1}-h_j)]\xi_j=0,
\end{eqnarray}
and
\begin{eqnarray}\label{42}
\Big[-\frac{c^3}{6}(h_{j+1}-h_j)\Big]\zeta_j+\Big[-\frac{\tau c^3}{2}\Big]\xi_j=0.
\end{eqnarray}
Multiplying both equations by $\frac{2}{c^3}$ gives
\begin{eqnarray}\label{43}
\tau \zeta_j+ \frac{1}{3}(h_{j+1}-h_j)\xi_j=0,
\end{eqnarray}
and
\begin{eqnarray}\label{44}
\frac{1}{3}(h_{j+1}-h_j)\zeta_j+\tau \xi_j=0.
\end{eqnarray}
Equations (\ref{43}) and (\ref{44}) can be written in a matrix system as
\begin{equation}\label{45}
\left[
\begin{array}{cc}
\tau & \frac{1}{3}(h_{j+1}-h_j) \\
\frac{1}{3}(h_{j+1}-h_j) & \tau
\end{array}
\right]
 \left[
\begin{array}{c}
\zeta_j \\
\xi_j
\end{array}
\right]
=
  \left[
\begin{array}{c}
0 \\
0
\end{array}
\right]\, .
\end{equation}
Note that since not all $\zeta_j$ and $\xi_j$ are zeros for all $j$, then it is clear that
\begin{equation}\label{46}
det\left[
\begin{array}{cc}
\tau & \frac{1}{3}(h_{j+1}-h_j) \\
\frac{1}{3}(h_{j+1}-h_j) & \tau
\end{array}
\right]=0,
\end{equation}
where $det(D)$ is the determinant of the matrix $D$. Solving Equation (\ref{46}) leads to
\begin{equation}\label{47}
\tau_j=\frac{1}{3}(h_{j+1}-h_j) .
\end{equation}
which is the desired result.\hfill{$\blacksquare$}\\

\section{Numerical Results and Discussion}
To make the discussion more beneficial and clearer, and in the spirit of fair comparison, we will compare the results of the stability scheme presented here to the computational results of \cite{ALMAN,ALMANA,ALMA}. The computation is carried out for the Hydrogen-like Ununoctium ion where the atomic number and atomic weight are respectively 118 and 294. The computation is majorally performed for the point nucleus for which the approximated eigenvalues can be compared with the exact eigenvalues obtained by the relativistic formula (\ref{7}). For the case of extended nucleus, we will, as mentioned before, assume uniformly distributed charge in the region $[0, R]$, where $R$ is the radius of the nucleus. The intensity of the nodes distribution near the nucleus is controlled by the parameter $\eps$ that plays a major role in Formula (\ref{8}), where the most appropriate values of $\eps$ are those that are living in $[10^{-6}\,,\,10^{-4}]$, see \cite{ALMAN}. For all of the computational results below, we have assumed that the nodes intensity parameter  $\eps=10^{-4}$.\\

Below, the computation is considered with $\kappa=\pm2, \pm3,\ldots$, the general case where the spinors are vanishing smoothly with zero derivatives at the boundaries, i.e., homogeneous Dirichlet and homogeneous Neumann boundary conditions. The case when $\kappa=\pm1$, and as mentioned before, there is no differences in the computation but nonhomogeneous Neumann boundary condition should instead be considered. Thus, in the computation of these two cases, small modification should be considered near the boundaries in the programming code. However, for a general discussion, we will consider the general case, that is when $\kappa=\pm2, \pm3,\ldots$.\\

\begin{table}[h]
\begin{footnotesize}
\caption{The first computed eigenvalues of the electron in the Hydrogen-like Ununoctium ion for $\kappa=-2$ using the usual FEM and the stability scheme with linear basis function for point nucleus, where the number of nodes is $600$, and $\eps=10^{-4}$. }
\centering
\begin{tabular}{@{} l c c r @{} }
\hline\hline
Level & Usual FEM & Stabilized FEM & Exact solution \\
& Linear basis & Linear basis & Relativistic Formula \\ [0.5ex]
\hline
1&-1829.630750908&-1829.630678009&-1829.630750908\\
2&-826.7683699234&-826.7681327991&-826.7683539069\\
3&-463.1183759679&-463.1178925700&-463.1183252634\\
&\cellcolor[gray]{0.6}-294.6216782193&&Spurios Eigenvalue\\
4&-294.4510822666&-294.4502765309&-294.4509801141\\
5&-203.2421234746&-203.2409198509&-203.2419549027\\
6&-148.5536893591&-148.5520121218&-148.5534402360\\
&\cellcolor[gray]{0.6}-113.4611501523&&Spurios Eigenvalue\\
7&-113.2482614926&-113.2460345755&-113.2479180697\\
8&-89.15839677745&-89.15554369439&-89.15794547564\\
9&-71.99903774457&-71.99548153219&-71.99846504808\\
&\cellcolor[gray]{0.6}-59.57649074983&&Spurios Eigenvalue\\
10&-59.34933184120&-59.34499500331&-59.34862423729\\
11&-49.75886521413&-49.75366967052&-49.75800915710\\
12&-42.31613542297&-42.31000245862&-42.31511730902\\
&\cellcolor[gray]{0.6}-36.65876644972&&Spurios Eigenvalue\\
13&-36.42517755182&-36.41802776855&-36.42398370073\\
14&-31.68311361412&-31.67486688495&-31.68173025393\\
15&-27.80972133834&-27.80029676250&-27.80813459180\\
[1ex]
\hline\hline
\end{tabular}
\label{table:nonlin}
\end{footnotesize}
\end{table}
In Tables 2 and 3, the usual and the stabilized FEM with linear basis functions are applied for approximating the eigenvalues of the radial Dirac operator with $\kappa=-2$ and $\kappa=2$ respectively for point nucleus. The number of nodes used is $600$ and the computation is carried out for the Hydrogen-like Ununoctium ion. The so called instilled spurious eigenvalues (the gray-colored ones except the first value on the top of the second column of Table 3) clearly presented in the computation using the usual FEM. Also the spurious eigenvalues caused by the unphysical coincidence phenomenon (the gray-colored value on the top of the second column of Table 3) clearly presented in the computation using the usual FEM. While the computation of the eigenvalues using the stable finite element scheme is cleaned from both categories of the spectrum pollution.
\begin{table}[h]
\begin{footnotesize}
\caption{The first computed eigenvalues of the electron in the Hydrogen-like Ununoctium ion for $\kappa=2$ using the usual FEM and the stability scheme with linear basis function for point nucleus, where the number of nodes is $600$ and $\eps=10^{-4}$.}
\centering
\begin{tabular}{@{} l c c r @{} }
\hline\hline
Level & Usual FEM & Stabilized FEM & Exact solution \\
& Linear basis & Linear basis & Relativistic Formula \\ [0.5ex]
\hline
&\cellcolor[gray]{0.6}-1829.630750908&&Spurios Eigenvalue\\
1&-826.7683699236&-826.7682977877&-826.7683539068\\
2&-463.1183759680&-463.1181147468&-463.1183252633\\
&\cellcolor[gray]{0.6}-294.6216782190&   &Spurios Eigenvalue\\
3&-294.4510822666&-294.4505263188&-294.4509801141\\
4&-203.2421234749&-203.2411849302&-203.2419549026\\
5&-148.5536893591&-148.5522865624&-148.5534402360\\
&\cellcolor[gray]{0.6}-113.4611501522&   &Spurios Eigenvalue\\
6&-113.2482614926&-113.2463152035&-113.2479180697\\
7&-89.15839677744&-89.15582867204&-89.15794547563\\
8&-71.99903774456&-71.99576972929&-71.99846504808\\
&\cellcolor[gray]{0.6}-59.57649074972&   &Spurios Eigenvalue\\
9&-59.34933184115&-59.34528569217&-59.34862423728\\
10&-49.75886521408&-49.75396236511&-49.75800915710\\
11&-42.31613542286&-42.31029682491&-42.31511730902\\
&\cellcolor[gray]{0.6}-36.65876644981&   &Spurios Eigenvalue\\
12&-36.42517755184&-36.41832357217&-36.42398370072\\
13&-31.68311361413&-31.67516395745&-31.68173025392\\
14&-27.80972133816&-27.80059498195&-27.80813459179\\
[1ex]
\hline\hline
\end{tabular}
\label{table:nonlin}
\end{footnotesize}
\end{table}

In Table 4, a comparison between the $hp$-cloud Petrov-Galerkin ($hp$-CPG) method \cite{ALMAN} and the stabilized FEM with Linear basis functions is considered. The computation is carried out for the Hydrogen-like Ununoctium ion for point nucleus for $\kappa=-2$ and the number of nodal points used is $600$. The $hp$-CPG is obtained at $\rho_j= 2.2h_{j+1}$,  where the clouds are enriched by $P^t(x)=[1\,,\,x(1-x/2) \,\exp(-x/2)]$, see \cite{ALMAN}.
\begin{table}[h]
\begin{footnotesize}
\caption{The first computed eigenvalues of the electron in the Hydrogen-like Ununoctium ion for $\kappa=-2$ using the $hp$-cloud Petrov-Galerkin stability scheme (as in \cite{ALMAN}) and the current stability scheme using linear basis functions for point nucleus, the number of nodes is $600$.}
\centering
\begin{tabular}{@{} l c c c c r @{} }
\hline\hline
Level & $hp$-cloud & FEM & Exact solution &  Relative Error & Relative Error\\
& Petrov-Galerkin & Linear (Hat) & Relativistic Formula & hp-CPG & FEM\hspace{0.5cm} \\ [0.5ex]
\hline
1&-1829.628962027&-1829.630678009&-1829.630750908&0.0000009777&0.0000000398\\
2&-826.7707399381&-826.7681327991&-826.7683539069&0.0000028859&0.0000002674\\
3&-463.1232256402&-463.1178925700&-463.1183252634&0.0000105812&0.0000009343\\
4&-294.4572672676&-294.4502765309&-294.4509801141&0.0000213521&0.0000023894\\
5&-203.2490442643&-203.2409198509&-203.2419549027&0.0000348813&0.0000050927\\
6&-148.5610131628&-148.5520121218&-148.5534402360&0.0000509777&0.0000096134\\
7&-113.2557872543&-113.2460345755&-113.2479180697&0.0000694863&0.0000166316\\
8&-89.16599265476&-89.15554369439&-89.15794547564&0.0000902575&0.0000269384\\
9&-72.00661059443&-71.99548153219&-71.99846504808&0.0001131350&0.0000414386\\
10&-59.35681135728&-59.34499500331&-59.34862423729&0.0001379496&0.0000611511\\
11&-49.76619519073&-49.75366967052&-49.75800915710&0.0001645169&0.0000872118\\
12&-42.32326870782&-42.31000245862&-42.31511730902&0.0001926356&0.0001208752\\
13&-36.43207300156&-36.41802776855&-36.42398370073&0.0002220872&0.0001635167\\
14&-31.68973420172&-31.67486688495&-31.68173025393&0.0002526360&0.0002166349\\
15&-27.81603295481&-27.80029676250&-27.80813459180&0.0002840306&0.0002818538\\
[1ex]
\hline\hline
\end{tabular}
\label{table:nonlin}
\end{footnotesize}
\end{table}

\vspace{0.5cm}
In Table 5, a comparison between the stabilized FEM using the cubic hermitian \cite{ALMANA,ALMA} and using the linear basis functions is presented. The computation is carried out for the Hydrogen-like Ununoctium ion for point nucleus for $\kappa=-2$ and the number of nodal points used is $600$.

\begin{table}[h]
\begin{footnotesize}
\caption{The first computed eigenvalues of the electron in the Hydrogen-like Ununoctium ion for $\kappa=-2$ using the finite element stability scheme with cubic hermitian basis functions (as in \cite{ALMANA,ALMA}) and the current stability scheme using linear basis functions for point nucleus, the number of nodes is $600$.}
\centering
\begin{tabular}{@{} l c c c c r @{} }
\hline\hline
Level & FEM & FEM & Exact solution &  Relative Error & Relative Error\\
& Hermitian Cubic & Linear (Hat) & Relativistic Formula & Hermitian Cubic  & Linear (Hat) \\ [0.5ex]
\hline
1&-1829.630750699&-1829.630678009&-1829.630750908&0.0000000001142&0.0000000398\\
2&-826.7683538119&-826.7681327991&-826.7683539069&0.0000000001149&0.0000002674\\
3&-463.1183252175&-463.1178925700&-463.1183252634&0.0000000000991&0.0000009343\\
4&-294.4509800935&-294.4502765309&-294.4509801141&0.0000000000699&0.0000023894\\
5&-203.2419548930&-203.2409198509&-203.2419549027&0.0000000000477&0.0000050927\\
6&-148.5534402320&-148.5520121218&-148.5534402360&0.0000000000269&0.0000096134\\
7&-113.2479180654&-113.2460345755&-113.2479180697&0.0000000000379&0.0000166316\\
8&-89.15794546761&-89.15554369439&-89.15794547564&0.0000000000900&0.0000269384\\
9&-71.99846503277&-71.99548153219&-71.99846504808&0.0000000002126&0.0000414386\\
10&-59.34862421008&-59.34499500331&-59.34862423729&0.0000000004584&0.0000611511\\
11&-49.75800911278&-49.75366967052&-49.75800915710&0.0000000008907&0.0000872118\\
12&-42.31511724130&-42.31000245862&-42.31511730902&0.0000000016003&0.0001208752\\
13&-36.42398360216&-36.41802776855&-36.42398370073&0.0000000027061&0.0001635167\\
14&-31.68173011572&-31.67486688495&-31.68173025393&0.0000000043624&0.0002166349\\
15&-27.80813440400&-27.80029676250&-27.80813459180&0.0000000067534&0.0002818538\\
[1ex]
\hline\hline
\end{tabular}
\label{table:nonlin}
\end{footnotesize}
\end{table}

\vspace{1cm}
In Figure 1, the convergence rates for the computation of the first fifteen eigenvalues using the three methods, $hp$-CPG, stabilized FEM with linear basis, and stabilized FEM with cubic hermitian, are shown. It is clearly noticed that the convergence rate of the approximation using the FEM with cubic hermitian is better than those of the other two methods, while the approximation using the FEM with linear basis functions is better than the convergence rate of the $hp$-CPG method.

\begin{figure}[ht]
\centering
\includegraphics[width=12cm]{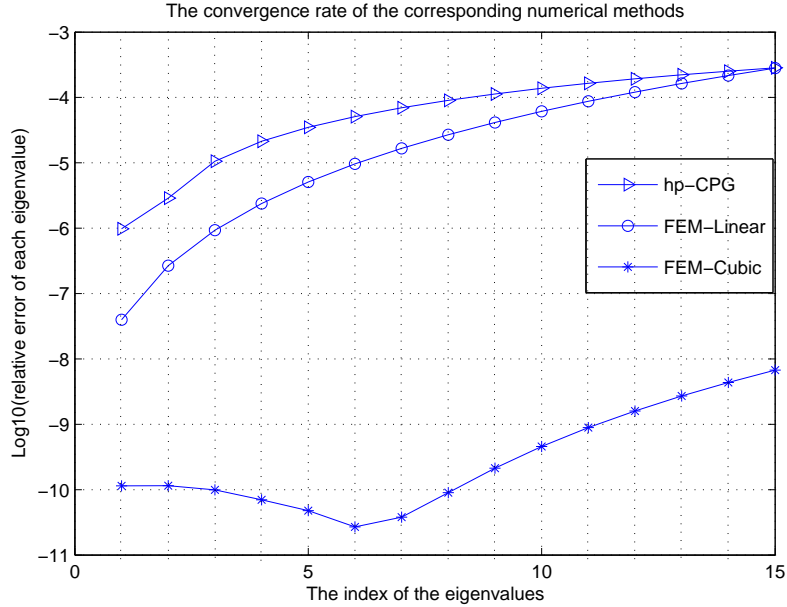}
\caption{Comparison between the the previously derived stability schemes (as in \cite{ALMAN,ALMANA,ALMA}) and the current stability scheme with linear basis functions, see the corresponding tables (Tables 4 and 5).}
\end{figure}

Table 6 presents the computation of the eigenvalues of the Hydrogen-like Ununoctium ion for $\kappa=-2$ for point nucleus with different numbers of nodal points.

\begin{table}[h]
\begin{footnotesize}
\caption{The first computed eigenvalues of the electron in the Hydrogen-like Ununoctium ion for $\kappa=-2$ for point nucleus using different numbers of nodes.}
\centering
\begin{tabular}{@{} l c c c c c r @{} }
\hline\hline
Level & $n=200$ & $n=400$ & $n=600$ &  $n=800$ &  $n=1000$ & Exact solution \\ [0.5ex]
\hline
1 &-1829.624974&-1829.630384&-1829.630678&-1829.630727&-1829.630741&-1829.630750 \\
2 &-826.7507746&-826.7672405&-826.7681327&-826.7682837&-826.7683250&-826.7683539 \\
3 &-463.0838205&-463.1161451&-463.1178925&-463.1181879&-463.1182689&-463.1183252 \\
4 &-294.3946426&-294.4474328&-294.4502765&-294.4507569&-294.4508885&-294.4509801 \\
5 &-203.1586471&-203.2367320&-203.2409198&-203.2416267&-203.2418202&-203.2419549 \\
6 &-148.4377973&-148.5462267&-148.5520121&-148.5529875&-148.5532545&-148.5534402 \\
7 &-113.0943605&-113.2383931&-113.2460345&-113.2473213&-113.2476733&-113.2479180 \\
8 &-88.96068950&-89.14578334&-89.15554369&-89.15718489&-89.15763368&-89.15794547 \\
9 &-71.75154013&-71.98333473&-71.99548153&-71.99752076&-71.99807804&-71.99846504 \\
10&-59.04590551&-59.33018952&-59.34499500&-59.34747626&-59.34815388&-59.34862423 \\
11&-49.39327045&-49.73592836&-49.75366967&-49.75663740&-49.75744729&-49.75800915 \\
12&-41.88210950&-42.28904304&-42.31000245&-42.31350159&-42.31445573&-42.31511730 \\
13&-35.91654123&-36.39356270&-36.41802776&-36.42210370&-36.42321419&-36.42398370 \\
14&-31.09390617&-31.64660327&-31.67486688&-31.67956554&-31.68084455&-31.68173025 \\
15&-27.13436106&-27.76793634&-27.80029676&-27.80566461&-27.80712439&-27.80813459 \\
[1ex]
\hline\hline
\end{tabular}
\label{table:nonlin}
\end{footnotesize}
\end{table}

Figure 2 shows the convergence rate of the approximation for the first five eigenvalues that are presented in Table 6.

\begin{figure}[h]
\centering
\includegraphics[width=14cm]{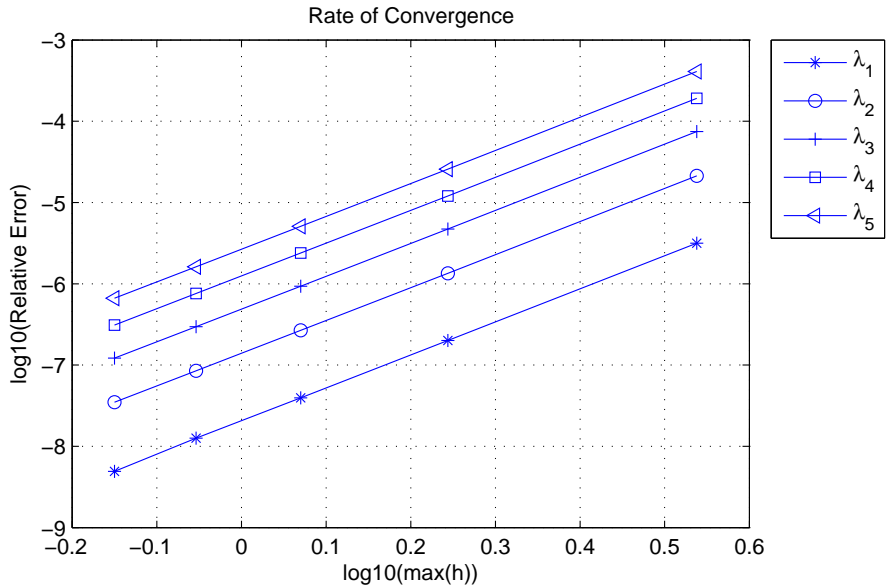}
\caption{The convergence rate of the approximation of the first five eigenvalues of the electron in the Ununoctium ion (Table 6)using the stable FEM with linear basis functions.}
\end{figure}

In Table 7 below, the computation is carried out for extended nucleus with different values of $\kappa$. The number of nodes used is 600, where 40 out of them is in the range $[0, R]$, and 560 nodal points in the rest of the domain.\\

\textbf{Conclusion.}
To conclude the work done of this work, we will summarize the subtle points as a comparison between the presented paper and the works done in \cite{ALMAN,ALMA}:

\begin{itemize}
\item As of the schemes presented in \cite{ALMAN,ALMA}, the new stability scheme provides a complete remedy of the spurious eigenvalues of both categories (instilled spurious eigenvalues and the spuriosity caused by the so-called unphysical coincidence phenomenon), for all Hydrogen-like ions, for all values of the quantum number $\kappa$, and for both point and extended nucleus.\\

\item The derivation of the stability parameter $\tau$ here is simpler, faster, and less time consuming compared with the derivations of $\tau$ in each of \cite{ALMAN} and \cite{ALMA}. Moreover, the derivation of $\tau$ does not require the fact of the accumulated eigenvalues of the radial Dirac operator.\\

\item The rate of convergence of the stabilized FEM using the cubic hermitian basis functions is much better than both stabilized FEM using the linear basis functions and the $hp$-CPG method. On the other hand, the stabilized finite element method using the linear basis functions is relatively better than the $hp$-CPG method.
\end{itemize}

\begin{table}[ht]
\begin{footnotesize}
\caption{The first computed eigenvalues of the electron in the Hydrogen-like Ununoctium ion using the stability scheme for extended nucleus with different values of $\kappa$.}
\centering
\begin{tabular}{@{} l c c c c c c c r @{} }
\hline\hline
Level & $\kappa=-2$ & $\kappa=2$ & $\kappa=-3$ &  $\kappa=3$ & $\kappa=-4$ & $\kappa=4$ & $\kappa=-5$ & $\kappa=5$\\[0.5ex]
\hline\hline
1 &-1829.6307&		    &	    	&	    	&		    &		    &		    &		    \\
2 &-826.76812&-826.76830&		    &	     	&		    &		    &		    &	    	\\
3 &-463.11788&-463.11811&-790.18014	&		    &		    &		    &		    &		    \\
4 &-294.45026&-294.45051&-447.43111	&-447.43131	&		    &		    &		    &		    \\
5 &-203.24089&-203.24116&-286.40405	&-286.40434	&-440.28637	&		    &		    &		    \\
6 &-148.55197&-148.55226&-198.59429	&-198.59463	&-282.71295	&-282.71318	&		    &		    \\
7 &-113.24598&-113.24627&-145.63495	&-145.63533	&-196.45227	&-196.45261	&-280.57597	&		    \\
8 &-89.155479&-89.155772&-111.29810	&-111.29850	&-144.28551	&-144.28592	&-195.20917	&-195.20940	\\
9 &-71.995402&-71.995698&-87.791571	&-87.791981	&-110.39463	&-110.39509	&-143.50105	&-143.50143	\\
10&-59.344898&-59.345196&-71.003860	&-71.004280	&-87.157686	&-87.158179	&-109.86879	&-109.86926	\\
11&-49.753553&-49.753854&-58.601791	&-58.602219	&-70.542312	&-70.542830	&-86.788430	&-86.788962	\\
12&-42.309865&-42.310168&-49.182450	&-49.182885	&-58.255460	&-58.255996	&-70.273285	&-70.273862	\\
13&-36.417868&-36.418172&-41.861579	&-41.862019	&-48.916026	&-48.916575	&-58.053510	&-58.054121	\\
14&-31.674683&-31.674988&-36.059596	&-36.060039	&-41.652292	&-41.652852	&-48.760634	&-48.761272	\\
15&-27.800086&-27.800392&-31.383873	&-31.384320	&-35.892236	&-35.892805	&-41.530218	&-41.530877	\\
16&-24.594310&-24.594617&-27.560825	&-27.561275	&-31.247969	&-31.248546	&-35.794625	&-35.795300	\\
17&-21.911846&-21.912154&-24.395109	&-24.395562	&-27.448980	&-27.449563	&-31.168722	&-31.169411	\\
18&-19.644685&-19.644995&-21.744246	&-21.744702	&-24.301976	&-24.302564	&-27.383783	&-27.384484	\\
19&-17.711297&-17.711607&-19.502347	&-19.502804	&-21.665885	&-21.666478	&-24.247713	&-24.248424	\\
20&-16.049205&-16.049517&-17.589391	&-17.589850	&-19.435800	&-19.436397	&-21.620256	&-21.620976	\\
21&-14.609890&-14.610202&-15.944004	&-15.944465	&-17.532407	&-17.533008	&-19.397080	&-19.397807	\\
22&-13.355203&-13.355516&-14.518478	&-14.518941	&-15.894842	&-15.895447	&-17.499280	&-17.500014	\\
23&-12.254823&-12.255137&-13.275274	&-13.275739	&-14.475776	&-14.476384	&-15.866293	&-15.867032	\\
24&-11.284401&-11.284715&-12.184532	&-12.184999	&-13.237954	&-13.238565	&-14.451007	&-14.451752	\\
25&-10.424204&-10.424520&-11.222260	&-11.222728	&-12.151732	&-12.152346	&-13.216336	&-13.217086	\\
26&-9.6581053&-9.6584214&-10.369002	&-10.369471	&-11.193283	&-11.193900	&-12.132761	&-12.133516	\\
27&-8.9728188&-8.9731357&-9.6088459	&-9.6093166	&-10.343281	&-10.343900	&-11.176552	&-11.177311	\\
28&-8.3573215&-8.3576390&-8.9286794	&-8.9291515	&-9.5859147	&-9.5865361	&-10.328457	&-10.329220	\\
29&-7.8024051&-7.8027234&-8.3176166	&-8.3180900	&-8.9081527	&-8.9087763	&-9.5727260	&-9.5734926	\\
30&-7.3003298&-7.3006487&-7.7665602	&-7.7670349	&-8.2991732	&-8.2997989	&-8.8963737	&-8.8971438	\\
31&-6.8445513&-6.8448708&-7.2678602	&-7.2683361	&-7.7499308	&-7.7505585	&-8.2886159	&-8.2893893	\\
32&-6.4295064&-6.4298266&-6.8150462	&-6.8155231	&-7.2528176	&-7.2534472	&-7.7404377	&-7.7412140	\\
33&-6.0504418&-6.0507625&-6.4026154	&-6.4030934	&-6.8013976	&-6.8020290	&-7.2442557	&-7.2450349	\\
34&-5.7032762&-5.7035974&-6.0258647	&-6.0263437	&-6.3901966	&-6.3908296	&-6.7936540	&-6.7944358	\\
35&-5.3844899&-5.3848116&-5.6807547	&-5.6812345	&-6.0145348	&-6.0151693	&-6.3831751	&-6.3839594	\\
36&-5.0910346&-5.0913567&-5.3638007	&-5.3642813	&-5.6703924	&-5.6710283	&-6.0081529	&-6.0089395	\\
37&-4.8202596&-4.8205820&-5.0719839	&-5.0724651	&-5.3543012	&-5.3549383	&-5.6645791	&-5.6653678	\\
38&-4.5698511&-4.5701737&-4.8026783	&-4.8031601	&-5.0632560	&-5.0638943	&-5.3489950	&-5.3497856	\\
39&-4.3377824&-4.3381052&-4.5535914	&-4.5540737	&-4.7946428	&-4.7952820	&-5.0584037	&-5.0591959	\\
40&-4.1222722&-4.1225951&-4.3227147	&-4.3231974	&-4.5461789	&-4.5468189	&-4.7901978	&-4.7909915	\\ [1ex]
\hline\hline
\end{tabular}
\label{table:nonlin}
\end{footnotesize}
\end{table}

\end{document}